\baselineskip=14pt
\parskip=10pt

\magnification=\magstephalf

\def\1{{\overline{1}}}
\def\2{{\overline{2}}}
\parindent=0pt
\overfullrule=0in

\def\frac#1#2{{#1 \over #2}}
\bf
\centerline
{
The Reciprocal of $\sum _{n\geq 0}a^nb^n$  for non-commuting $a$ and $b$,
} 
\centerline 
{Catalan numbers
and non-commutative quadratic equations }
\rm
\bigskip
\centerline{ \it
Arkady BERENSTEIN, Vladimir RETAKH, Christophe REUTENAUER  and Doron ZEILBERGER
}

{\bf Note:} This article is accompanied by the Maple package {\tt NCFPS} downloadable from \hfill\break
{\tt http://www.math.rutgers.edu/\~{}zeilberg/tokhniot/NCFPS}

The aim of this paper is to describe the inversion of the sum $\sum _{n\geq 0}a^nb^n$ where $a$ and $b$
are non-commuting variables as a formal series in $a$ and $b$. We show that the inversion satisfies
a non-commutative quadratic equation and that the number of certain monomials in its homogeneous components
equals a Catalan number. We also study general solutions of similar quadratic equations. 

\smallskip
\centerline {\bf 1. Inverting $\sum _{n\geq 0}a^nb^n$.}
\bigskip
Our goal is to find an inverse of the series $\sum _{n\geq 0}a^nb^n$ where $a$ and $b$ are non-commuting
variables. The answer to this question is given by the following theorem.

Let $a,b,x$ be (completely!) {\bf non-commuting} variables (``indeterminates"). Define
a sequence of polynomials $d_n(a,b,x)$ ( $n \geq 1$) recursively as follows:
$$
d_1(a,b,x)=1 \quad , \ \ (1a)
$$
$$
d_{n}(a,b,x)= d_{n-1}(a,b,x) x+ \sum_{k=2}^{n-1}  d_{n-k}(a,b,x)\, a \, d_k(a,b,x) \, b\quad  (n \geq 2) \quad . (1b)
$$
Also define the sequence of polynomials $c_n(a,b,x)$ as follows:
$$
c_n(a,b,x)=a\, d_n(a,b,x)\, b \quad (n \geq 1) \quad .
$$

{\bf Theorem 1}:
$$
1-\sum_{n=1}^{\infty} c_n(a,b,ab-ba) \, = \, \left ( \sum_{n \geq 0} a^nb^n \right )^{-1} \quad .
$$

It follows immediately that the number of monomials in $a,b$ and $x$ in the polynomial $d_n(a,b,x)$ 
is the $(n-1)$-th Catalan number. In particular, $d_1=1$, $d_2=x$, $d_3=x^2+axb$,
$$
d_4=x^3+ax^2b+axbx+xaxb+a^2xb^2,
$$
$$
d_5=x^4+ax^2bx+axbx^2+xaxbx+a^2xb^2x+x^2axb+axbaxb
$$
$$
+xax^2b+xaxb^2+ax^3b+a^2x^2b^2+a^2xbxb+axaxb^2+a^3xb^3.
$$

We will give an algebraic and a combinatorial proof of the theorem.
A simple algebraic proof is based on two lemmas.

{\bf Lemma 2}: Let $S$ be a formal series in $a$ and $b$ such that $S=1+aSb$. Observe that the inverse of $S$
is of the form $1-C$ where $C=aDb$ and the series $D$ satisfies the equation
$$D=1+D(x-ab)+DaDb \ \ \ \ (2)$$
and $x=ab-ba$.

{\bf Proof}: We are looking for the inverse of $S$ in the form
$1-C$ where $C=aDb$. 

We have
$$CS=(1-S^{-1})S=S-1=aSb.$$

Hence $$C(1+aSb)=aSb,$$
$$C+CaSb=aSb,$$
$$aDb+aDbaSb=aSb.$$

So,

$$D+DbaS=S$$ and $$D(1+baS)=S$$ or
$$D(S^{-1}+ba)=1.$$

It implies that $$D(1-C+ba)=1$$
and
$$D=1+DaDb-Dba$$ which immediately implies equation (2). 

{\bf Lemma 3}: Let the degree of indeterminates $a$ and $b$ in equation (2) equal one and the degree of $x$ equal two.
Then the solution of equation (2) is given by formula
$$D=\sum _{n\geq 1}d_n(a,b,x)$$
where polynomials $d_n(a,b,x)$ satisfy equations  (1).

{\bf Proof}: 
Note that $D=\sum _{n=1}^{\infty}d_n$ where $d_n=d_n(a,b,x)$ are homogeneous polynomials
in $a$ and $b$ of degree $2n-2$, $n=1,2,\dots $.

The terms of degree $0$ and $2$ are: $d_1=1$ and $d_2=x.$

Take the term of degree $2n-2$, $n\geq 3$:
$$d_n=d_{n-1}(x-ab)+\sum _{k=1}^{n-1}d_{n-k}ad_kb
=d_{n-1}(x-ab)+d_{n-1}ab + d_1ad_{n-1}b+\sum _{k=2}^{n-2}d_{n-k}c_k=$$
$$=d_{n-1}x+ad_{n-1}b+\sum _{k=2}^{n-2}d_{n-k}c_k.$$
QED

Let $S=\sum _{n\geq 0}a^nb^n$. Then $S$ satisfies equation $S=1+aSb$ and Theorem 1 follows from Lemmas 2 and 3.
\bigskip  
{\bf Combinatorial Proof}: Consider the set of {\it lattice walks} in the 2D rectangular lattice, starting
at the origin, $(0,0)$ and ending at $(n-1,n-1)$, where one can either make a {\it horizontal} step
$(i,j) \rightarrow (i+1,j)$,
(weight $a$), a {\it vertical} 
step $(i,j) \rightarrow (i,j+1)$, (weight $b$) or a diagonal step 
$(i,j) \rightarrow (i+1,j+1)$, (weight $x$), 
always staying in the region $i \geq j$, and
where you can never have a horizontal step followed immediately by a vertical step. In other words, you may never
venture to the region $i<j$, and you can never have the  Hebrew letter Nun (alias the mirror-image of the
Latin letter $L$) when you draw the path on the plane. The weight of a path is the
product (in order!) of the weights of the individual steps.

For example, when $n=2$ the only possible path is $(0,0) \rightarrow (1,1)$,
whose weight is $x$.

When $n=3$ we have two paths.
The path $(0,0) \rightarrow (1,1) \rightarrow (2,2)$ whose weight is $x^2$ and
the path $(0,0) \rightarrow (1,0) \rightarrow (2,1) \rightarrow (2,2)$ whose weight is $axb$.

When $n=4$ we have five paths:

The path $(0,0) \rightarrow (1,1) \rightarrow (2,2) \rightarrow (3,3) $ whose weight is $x^3$,

the path $(0,0) \rightarrow (1,0) \rightarrow (2,1) \rightarrow (3,2) \rightarrow (3,3) $ whose weight is $ax^2b$,

the path $(0,0) \rightarrow (1,0) \rightarrow (2,1) \rightarrow (2,2) \rightarrow (3,3) $ whose weight is $axbx$,

the path $(0,0) \rightarrow (1,1) \rightarrow (2,1) \rightarrow (3,2) \rightarrow (3,3) $ whose weight is $xaxb$, and

the path $(0,0) \rightarrow (1,0) \rightarrow (2,0) \rightarrow (3,1) \rightarrow (3,2) \rightarrow (3,3) $ whose weight is $a^2xb^2$.

It is very well-known, and rather easy to see, that the number of such paths are given by the Catalan numbers
$C(n-1)$, [2] {\tt http://oeis.org/A000108} .

We claim that the {\it weight-enumerator} of the set of such walks equals $d_{n}(a,b,x)$. Indeed, since
the walk ends on the diagonal, at the point $(n-1,n-1)$, the last step must be either a diagonal step

$$
(n-2,n-2) \rightarrow (n-1,n-1) \quad,
$$
whose weight-enumerator, by the inductive hypothesis is $d_{n-1}(a,b,x) x$, or else let $k$ be the smallest
integer such that the walk passed through $(n-k-1,n-k-1)$ (i.e. the penultimate encounter with the diagonal).
Note that $k$ can be anything between $2$ and $n-1$. The weight-enumerator of the set of paths from
$(0,0)$ to $(n-k-1,n-k-1)$ is $d_{n-k}(a,b,x)$, and the weight-enumerator of the set of paths from
$(n-k-1,n-k-1)$ to $(n-1,n-1)$ that never touch the diagonal, is $ad_{k}(a,b,x)b$. So the weight-enumerator is
$d_{n-k}(a,b,x)\, a \, d_k(a,b,x) b$ giving the above recurrence for $d_n(a,b,x)$.

It follows that $c_n(a,b,x)=a d_n(a,b,x) b$ is the weight-enumerator of all paths from
$(0,0)$ to $(n,n)$ as above with the additional property that except at the beginning
($(0,0))$ and the end ($(n,n)$) they always stay {\bf strictly} below the diagonal.

Now what does $c_n(a,b,ab-ba)$ weight-enumerate? Now there is a new rule in Manhattan, ``no shortcuts'',
one may not walk diagonally. So every diagonal step $(i,j) \rightarrow (i+1,j+1)$ must decide
whether 

to go first horizontally, and then vertically $(i,j) \rightarrow (i+1,j) \rightarrow (i+1,j+1)$, 
replacing $x$ by $ab$, or

to go first vertically, and then horizontally $(i,j) \rightarrow (i,j+1) \rightarrow (i+1,j+1)$, 
replacing $x$ by  $-ba$.

This has to be decided, independently for each of the diagonal steps that formerly had weight $x$.
So a path with $r$ diagonal steps gives rise to $2^r$ new paths with sign $(-1)^s$ where $s$ is the
number of places where it was decided to go through the second option.

So $c_n(a,b,ab-ba)$ is the weight-enumerators of pairs of paths
$[P,K]$
where $P$ is the original path featuring a certain number of diagonal steps $r$,
and $K$ is one of its $2^r$ ``children'', paths with only horizontal and vertical steps,
and weight $\pm weight(P)$, where we have a plus-sign if an even number of the $r$ diagonal
steps became {\it vertical-then-horizontal} (i.e. $ba$) and a minus-sign otherwise.

As we look at the weights of the children $K$, 
sometimes we have the same path coming from different
parents. Let's call a pair $[P,K]$ {\it bad} if the path $P$  has a ``$ba$'' {\it strictly-under}
the diagonal, i.e.  a ``vertical step followed by a horizontal step'' that does not touch the
diagonal. Write $K$ as $K=w_1(ba)^sw_2$ where $w_1$ does not have any sub-diagonal $ba$'s and $s$ is as large
as possible. Then the parent must be either of the form $P=W_1 x^s W_2$ where the $x^s$ corresponds to the $(ba)^s$,
or of the form  $P'=W_ 1 b x^{s-1} a W_2$. In the former case attach $[W_1 x^s W_2, K]$ to $[W_ 1 b x^{s-1} a W_2,K]$
and in the latter case  vice-versa. This is a weight-preserving and {\bf sign-reversing} involution among the
bad pairs, so they all kill each other.

It remains to weight-enumerate the {\it good pairs}. It is easy to see that the good pairs
are pairs $[P,K]$ where $K$ has the form $K=a^{i_1}b^{i_1}a^{i_2}b^{i_2} \dots a^{i_s}b^{i_s}$
for some $s \geq 1$ and integers $i_1, \dots, i_s \geq 1$ summing up to $n$
(this is called a {\it composition} of $n$). It is easy to see that for each such $K$, (coming from a good pair
$[P,K]$)there can only
be {\bf one} possible {\it parent} $P$. The sign of a good pair
$$
[P,a^{i_1}b^{i_1}a^{i_2}b^{i_2} \dots a^{i_s}b^{i_s} ] \quad ,
$$
is $(-1)^{s-1}$, since it touches the diagonal $s-1$ times, and each of these touching points
came from an $x$ that was turned into  $-ba$.

So $1-\sum_{n=1}^{\infty} c_n(a,b,ab-ba)$ turned out to be the sum of all the weights of compositions (vectors of positive integers)
$(i_1, \dots, i_s)$ with the weight $(-1)^{s} a^{i_1}b^{i_1} \cdots a^{i_s}b^{i_s}$ over {\it all} compositions, but the
same is true of
$$
\left ( \sum_{n \geq 0} a^nb^n \right )^{-1}=
\left ( 1+ \sum_{n \geq 1} a^nb^n \right )^{-1}=
1+\sum_{s=1}^{\infty} (-1)^{s} \left ( \sum_{n \geq 1} a^nb^n \right )^s \quad .
$$
QED!

\bigskip
{\bf 2. Inversion of $1-aDb$ in the general case}.
\bigskip
The following is a variant of of path's model used in Section 1.
Call {\it Dyck path} a path that starts at the origin,  ends on the $x$-axis,
that uses the steps $(1,1)$ (denoted by $a$) and $(1,-1)$ (denoted by $b$), and
that never goes below the $x$-axis. It is coded by a {\it Dyck word}, e.g. $aaababbabb$.
Formally, a Dyck word has as many $a$'s than $b$'s, and each prefix of it
has at least as many $a$'s as $b$'s.

If we replace, in each Dyck word, each occurrence of $ab$ by a letter $x$, and
sum all these words, then we obtain the series $D=\sum _{n\geq 1}d_n$
described in Section 1.

If we replace each $ab$ by a letter $x$, except those at level $0$, then we
obtain the series
$$1+aUb=1+\sum _{n\geq 1}au_nb.$$

For a series $Z$ set $Z^*:=(1-Z)^{-1}$. Then
$$(aDb)^*=1+aUb.$$

{\bf Theorem 4}. One has the equation
$$U=(1+aUb)(1+(x-ab+ba)U)$$
that completely defines $U$.

{\bf Proof}: We have $(1-aDb)^{-1}=1+aUb$,
thus $1-aDb=(-aUb)^*$.

The defining equation for $D$ is
$$D=1+(x-ab+aDb)D\ \ \ \ \ \ (3)$$
which is a symmetric version of equation (2); it follows from the Dyck path model, by writing $D=1+(x+a(D-1)b)D$. 
Let $D=1+d_1(a,b,x)+d_2(a,b,x)+\dots $ and polynomials $d_n(a,b,x)$
that satisfy equations (1) without any assumptions on $x$.

We have $$1-aDb=(-aUb)^*=1-aUb+(aUb)^2-(aUb)^3+\dots \quad .$$
Therefore,
$$aDb=aUb-(aUb)^2+(aUb)^3-\dots =a(U-UbaU+UbaUbaU-\dots)b$$
$$=aU(1-baU+(baU)^2-\dots )b$$
and
$$D=U(-baU)^*.$$

Note that (3) implies 
$$U(-baU)^*=1+(x-ab+aU(-baU)^*b)U(-baU)^*$$
therefore,
$$U=1+baU+(x-ab)U+aU(-baU)^*bU$$
$$=1+(x-ab+ba)U+aU(1-baU+baUbaU-\dots )bU$$
$$=1+(x-ab+ba)U+aUbU-aUbaUbU+aUbaUbaUbU-\dots$$
$$=1+(x-ab+ba)U+(-aUb)^*aUbU.$$

Hence
$$(1+aUb)U=1+aUb+(1+aUb)(x-ab+ba)U+aUbU$$
and 
$$U=1+aUb+(1+aUb)(x-ab+ba)U$$
$$=(1+aUb)(1+(x-ab+ba)U).$$

QED
\medskip
{\bf Remark 5}. If we put $x=ab-ba$ in the last equation,
then $U=1+aUb$ which implies $U=\sum _{n\geq 1}a^{n-1}b^{n-1}$ and
$1+aUb=\sum _{n\geq 0}a^nb^n$.

Note that Theorem 4 does not imply that all coefficients in $U$
as series in $a$, $b$ and $x$ are positive. However, simple computations
show that the inversion of the series
$1-aDb$ is written in the form $$1+au_1b+au_2b+...$$ where the degree of $u_n$ is
$2n-2$, $n\geq 1$ and
$$u_1=1,$$
$$u_2=ba + x,$$
$$u_3=(ba)^2+xba+bax+axb+x^2,$$
$$u_4=(ba)^3+x(ba)^2+baxba+(ba)^2x+
a^2xb+axb^2a+ba^2xb$$
$$+x^2ba+xbax+bax^2+ax^2b^2+axbx+xaxb+x^3,$$
and so on. The positivity follows from the path interpretation at the beginning of the section.

{\bf Problem 6}: How to write a recurrence relations on $u_n$ similar to
relations (1). It must imply that the number of terms for $u_n$ is the $n$-th Catalan number.
It also must show that if $x=ab-ba$ then $u_n=a^{n-1}b^{n-1}$.

We may set $x=1$ and get
$$u_1=1,\ \ u_2=ba+1,\ \ u_3=(ba)^2+2ba+ab+1,$$
$$u_4=(ba)^3+3(ba)^2+ab^2+ba^2b+a^2b^2+3ba+3ab+1.$$

{\bf Problem 7}: How to describe polynomials $u_n$ for this and other specializations?
Any relations with known polynomials? 

\bigskip
{\bf 3. The Quasideterminant of a Jacobi Matrix}

\bigskip
In this section we discuss solutions of noncommutative quadratic equation (2)
using quasideterminants. Recall ([1]) that quasideterminant
$|A|_{pq}$ of the matrix $A=(a_{ij})$, $i,j=1,2,\dots $ is defined
as follows. Let $A^{pq}$ be the submatrix of $A$ obtained from $A$
by removing its $p$-th row and $q$-th column. Denote by $r_p$ and $c_q$
be the $p$-th row and the $q$-th column of $A$ with element $a_{pq}$
removed. Assume that matrix $A^{pq}$ is invertible. Then
$$|A|_{pq}:=a_{pq}-r_p(A^{pq})^{-1}c_q.$$

Let now $A=(a_{ij})$, $i,j\geq 1$ be a Jacobi matrix, i.e.
$a_{ij}=0$ if $|i-j|>1$. Set $T=I-A$, where $I$ is the infinite identity matrix.
Recall that
$$|T|_{11}^{-1}=1+\sum a_{1j_1}a_{j_1j_2}a_{j_2j_3}\dots a_{j_k1}$$
where the sum is taken over all tuples $(j_1,j_2,\dots , j_k)$, $j_1,j_2,\dots , j_k\geq 1$, $k\geq 1$.

Also,
$$|T|_{11}=1-a_{11}-\sum a_{1j_1}a_{j_1j_2}a_{j_2j_3}\dots a_{j_k1}$$
where the sum is taken over all tuples $(j_1,j_2,\dots , j_k)$, $j_1,j_2,\dots , j_k>1$, $k\geq 1$.

Assume that the degree of all diagonal elements $a_{ii}$
is two and the degree
of all elements $a_{ij}$ such that $i\neq j$ is one. Then
$$|T|_{11}^{-1}=1+\sum _{n\geq 1}t_n\ \ \ \ \ (3)$$
where $t_n$ is homogeneous polynomial of degree $2n$ in variables $a_{ij}$.

In particular,
$$
t_1=a_{11}+a_{12}a_{21},
$$
$$
t_2=a_{11}^2+a_{11}a_{12}a_{21}+a_{12}a_{21}a_{11}+a_{12}a_{22}a_{21}+
(a_{12}a_{21})^2+a_{12}a_{23}a_{32}a_{21}.$$

Note that each  monomial corresponds, in a one-to-one way, to a ``Schr\"oder walk'' 
[2] {\tt http://oeis.org/A006318}, hence:

{\bf Proposition 8}: The number of monomials of $t_n$ is the $n$-th Large 
Schr\"oder Number.

If we set $a_{11}=0$ we get walks  obviously counted by the ``little''  
Schr\"oder numbers
[2] {\tt http://oeis.org/A001003}, hence:

{\bf Proposition 9}: Set $a_{11}=0$. Then the number of monomials in each $t_n$
is $A001003[n]$.

Let now $a,x,b$ be formal variables, the degree of $a$ and $b$ is one and the degree of $x$ is two.
Set $a_{ii}=x-ab$, $a_{i,i+1}=a, a_{i+1,i}=b$ for all $i$. By the definition of
quasideterminants, we have
$$|T|_{11}=1-x+ab-a|T|_{11}^{-1}b.$$

Denote $|T|_{11}^{-1}$ by $D$. Then last equation can be written as
$$D^{-1}=1-x+ab - aDb$$
or
$$D=1+D(x-ab)+DaDb$$
which is exactly our equation (2).  
\bigskip
{\bf Acknowledgement}

We wish to thank Matthew Russell for his careful reading and numerous corrections
and improvements.

A. Berenstein was partially supported by  NSF grant DMS-1101507. 
V. Retakh was supported in part by NSA grant H98230-11-1-0136. 
D. Zeilberger was supported in part by NSF grant  DMS-0901226.

\bigskip
{\bf References}
\bigskip
[1] I. Gelfand, V. Retakh, Determinant of matrices over noncommutative rings,
Funct. Anal.Appl., {\bf 25}(1991), no 2, pp. 91-102

[2] Neil Sloane, {\it On-Line Encyclopedia of Integer Sequences}, {\tt http://www.oeis.org}.

\bigskip
\hrule
\bigskip
Arkady Berenstein, Department of Mathematics, University of Oregon, Eugene, OR 97403 {\tt arkadiy [at] math [dot] uoregon [dot] edu }
\smallskip
Vladimir Retakh, Department of Mathematics, Hill Center- Busch Campus, Rutgers University,
110 Frelinghuysen Rd, Piscataway, NJ 08854-8019 , {\tt vretakh [at] math [dot] rutgers [dot] edu}
\smallskip
Christophe Reutenauer, 
Laboratoire de Combinatoire et d'informatique MathC)matique (LaCIM) Universit\'e du Qu\'ebec \'a  Montr\'eal,
CP 8888, Succ. Centre-ville, Montr\'eal (Qu\'eec) H3C 3P8 
{\tt Reutenauer [dot] Christophe [at] uqam [dot] ca}
\smallskip
Doron Zeilberger, Department of Mathematics, Hill Center- Busch Campus, Rutgers University,
110 Frelinghuysen Rd, Piscataway, NJ 08854-8019 , {\tt zeilberg [at] math [dot] rutgers [dot] edu}
\bigskip
\hrule
\bigskip
Written: June 18, 2012. This version: Oct. 26, 2012.
  
\end